\input amstex
\documentstyle{amsppt}

\input latexpicobjs.tex
\input pictexwd.tex
\font\tiny=cmr5
\font\llarge=cmr42

\magnification=\magstep1 \NoRunningHeads
\NoBlackBoxes

\topmatter

\title
Mixing constructions with infinite invariant measure and spectral multiplicities
 \endtitle
\author  Alexandre~I.~Danilenko and Valery V. Ryzhikov
\endauthor
\abstract
We introduce high staircase infinite measure preserving transformations and prove that they are mixing under a restricted growth condition.
This is used to (i) realize each subset $M\subset\Bbb N\cup \{\infty\}$ as the set of essential values of the multiplicity function for  the Koopman operator of a mixing  ergodic  infinite measure preserving transformation,
(ii) construct mixing power weakly mixing infinite measure preserving transformations,
(iii) construct mixing Poissonian automorphisms with a simple spectrum, etc.
\endabstract

\address
Max Planck Institute for Mathematics, Vivatsgasse 7, D-53111 Bonn, Germany
\endaddress
\address
Permanent address:
 Institute for Low Temperature Physics
\& Engineering of National Academy of Sciences of Ukraine, 47 Lenin Ave.,
 Kharkov, 61164, UKRAINE
\endaddress
\email alexandre.danilenko\@gmail.com
\endemail

\address
 Department of Mechanics and Mathematics, Lomonosov
Moscow State University, GSP-1, Leninskie Gory, Moscow, 119991, Russian
Federation
\endaddress
\email vryzh\@mail.ru
\endemail

\endtopmatter
\document

\head 0. Introduction
\endhead

Let $T$ be an  invertible measure preserving transformation of a $\sigma$-finite measure space $(X,\goth B,\mu)$ and $\mu(X)=\infty$.
By $U_T$ we denote the Koopman unitary operator associated with $T$:
$$
U_Tf=f\circ T, \quad f\in L^2(X,\mu).
$$
The set of essential values of the spectral multiplicity function of $U_T$ will be denoted by $\Cal M(T)$.
We note that $\Cal M(T)$ is a non-empty subset of $\Bbb N\cup\{\infty\}$.
In our previous paper \cite{DaR} we showed  that for each subset $M\subset\Bbb N$, there exists an ergodic conservative $T$ such that $\Cal M(T)=M$.
The main purpose of the present paper is to sharpen (and generalize) this result. We will show that $T$ can be chosen {\it mixing} or, equivalently, of {\it zero type} (see \cite{Aa}, \cite{DaS2}).
This means that $U_T^n\to 0$ weakly as $n\to\infty$.

\proclaim{Theorem 0.1} Given any subset $M\subset\Bbb N\cup\{\infty\}$, there exists a mixing ergodic conservative infinite measure preserving transformation $T$ such that $\Cal M(T)=M$.
\endproclaim

We first consider the case where $M$ does not contain $\{\infty\}$.
The idea of the proof is to consider the Cartesian product $S\times T$,  where $S$ is  a rigid transformation with $\Cal M(S)=M$ (constructed in \cite{DaR}) and $T$  is a mixing transformation such that $\Cal M(S\times T)=\Cal M(S)$ and $S\times T$ is ergodic. We note that the Cartesian product of a mixing system with any other (even non-ergodic) system is mixing.
We construct $T$ as a certain   limit of a sequence of transformations $(T_n)_{n>0}$ satisfying $U_{T_n}^{H_k^{(n)}}\to\delta_n U_{T_n}$ weakly as $k\to\infty$ along a subsequence $H_k^{(n)}$ of rigidity for $S$, where  $\delta_n$ is a sequence of positive reals tending to $0$.
The transformations $T_n$ are not mixing of course.
However they have a {\it mixing part} which occupies more and more space as $n\to\infty$. In the limit it  {\it fills} the entire space.
We thus consider this construction   as some {\it forcing of mixing}.
It remains to find a model for the mixing parts.
In the case of probability preserving systems, staircase transformations \cite{Ad} played this role (see \cite{Ry1}, \cite{Ry2}, \cite{Ag}, \cite{Da4}).
Therefore it seems natural to use ``infinite'' staircase systems for our purposes. Hence our first step is to show that the infinite staircases are mixing. Unfortunately, the {\it restricted growth condition on the sequence of cuts} essentially used by Adams in \cite{Ad} is incompatible with the infiniteness of the invariant measure.
That is why instead of pure staircases we introduce so-called {\it high} staircases.
Geometrically this means that on each step of the inductive cut-and-stack construction we insert a layer of spacers between the tower and the staircase {\it roof}.
If the layers are {\it sufficiently thick}, the corresponding high staircase transformation preserves an {\it infinite} measure.
At the same time an analogue of the restricted growth condition can hold for high staircases (see Definition~2.3 below).
Modifying Adams' argument from \cite{Ad} we show the following theorem,
 which is of independent interest.

\proclaim{Theorem 0.2}
Under the restricted growth condition each infinite measure preserving high staircase transformation
is mixing.
\endproclaim

We then utilize the transformations from Theorem~0.2 to prove Theorem~0.1
in the aforementioned way.
The case where $M\ni\infty$ comes to the above by considering a product
$T\times B$, where $B$ is a Bernoulli shift and $T$ is a  transformation from the claim of Theorem~1 with $\Cal M(T)=M\setminus\{\infty\}$.
 This idea `works' if the maximal spectral type of $T$ is singular.
We show that such $T$ exists.

 As a byproduct we can solve an  open problem related to weak mixing for infinite measure preserving systems. A transformation $T$  is called {\it power weakly mixing} if for each finite integer sequence $(n_1,\dots,n_k)$, $n_i\ne 0$ for all $i=1,\dots,k$, the product
$T^{n_1}\times \cdots\times T^{n_k}$ is ergodic.
A number of rank-one power weakly mixing transformations  with {\it exotic} properties are known so far
\cite{AFS}, \cite{Da1}, \cite{Da2}, \cite{DaS1} (see also  surveys \cite{Da3}, \cite{DaS2} and reference therein). However all of them are either rigid or partially rigid. Hence a problem arises:
\roster
{\it is there a mixing power weakly mixing rank-one infinite measure preserving map?}
\endroster
We answer affirmatively by showing the following theorem.

\proclaim{Theorem 0.3}
There is a power weakly mixing high staircase transformation satisfying the restricted growth condition.
\endproclaim

We  conjecture that every high staircase transformation is power weakly mixing.

The method of forcing  mixing  originated in \cite{Ry1} was used in \cite{Ry2} and  \cite{Ag} to construct  a  mixing rank-one finite measure preserving transformations $T$ such that the unitary
$$
\exp(\widehat U_T):=\bigoplus_{n=0}^\infty  \widehat U_T^{\odot n}
$$
has a simple spectrum.
Here $\widehat U_T$ denotes the restriction of $U_T$ onto the orthocomplement to the constant functions.
We  establish an infinite version of  this result.

\proclaim{Theorem 0.4} There is a mixing  rank-one conservative infinite measure preserving transformation $T$ such that the unitary operator $\exp(U_T)$ has a simple spectrum.
\endproclaim

This theorem has some applications to the theory of Poissonian automorphisms (see \cite{CFS}, \cite{Ne}, \cite{Ro1}, \cite{Ro2}).

\proclaim{Corollary 0.5}
There is a  mixing (of all orders) Poisson suspension with a simple spectrum.
\endproclaim

If a Poissonian automorphism does not have a simple spectrum then the set of its spectral multiplicities is infinite.
We  construct examples of mixing Poisson automorphisms $S$ such that  $\Cal M(S)=\{p,p^2,p^3,\dots\}$ (for an arbitrary $p>1$), and  also $\Cal M(S)=\{1, 3, 3\cdot 5,3\cdot 5\cdot 7,\dots\}$ (for other $S$, of course).

\head 1. Cut-and-stack and $(C,F)$-constructions
\endhead

To prove Main Theorem we will use  the $(C,F)$-construction (see \cite{dJ}, \cite{Da1}, \cite{Da3}).  We now briefly outline its formalism.
Let two sequences $(C_n)_{n>0}$ and $(F_n)_{n\ge 0}$ of finite subsets in $\Bbb Z$ are given such that:
\roster
\item"---"
$F_n=\{0,1,\dots,h_n-1\}$,  $\# C_n>1$, $0\in C_n$,
\item"---" $F_n+C_{n+1}\subset F_{n+1}$,
\item"---" $(F_{n}+c)\cap (F_n+c')=\emptyset$ if $c\ne c'$, $c,c'\in C_{n+1}$,
\item"---" $\lim_{n\to\infty}\frac{h_n}{\#C_1\cdots\# C_n}=\infty$.
\endroster
Let $X_n:=F_n\times C_{n+1}\times C_{n+2}\times\cdots$. Endow this set with the (compact Polish) product topology. The  following map
$$
(f_n,c_{n+1},c_{n+2},\dots)\mapsto(f_n+c_{n+1},c_{n+2},\dots)
$$
is a topological embedding of $X_n$ into $X_{n+1}$. We now set $X:=\bigcup_{n\ge 0} X_n$ and endow it with the (locally compact Polish) inductive limit topology. Given $A\subset F_n$, we denote by $[A]_n$ the following cylinder: $\{x=(f,c_{n+1},\dots,)\in X_n\mid f\in A\}$. Then $\{[A]_n\mid A\subset F_n, n>0\}$ is the family of all compact open subsets in $X$. It forms a base of the topology on $X$.

Let  $\Cal R$ stand for the {\it tail} equivalence relation on $X$: two points $x,x'\in X$ are $\Cal R$-equivalent if there is $n>0$ such that $x=(f_n,c_{n+1},\dots),\ x'=(f_n',c_{n+1}',\dots)\in X_n$ and $c_m=c_m'$ for all $m>n$. There is a   non-atomic Borel infinite $\sigma$-finite measure $\mu$ on $X$ which is invariant (and ergodic) under $\Cal R$ and such that $\mu(X_0)=h_0$.
It is unique up to scaling.

Now we define a transformation $T$ of $(X,\mu)$ by setting
$$
T(f_n,c_{n+1},\dots):=(1+f_n,c_{n+1},\dots )\text{ whenever }f_n<h_n-1,\ n>0.
$$
This formula defines $T$ partly on $X_n$. When $n\to\infty$, $T$ extends to the entire $X$ (minus countably many points) as a $\mu$-preserving invertible transformation. Moreover, the $T$-orbit equivalence relation  coincides  with $\Cal R$ (on the subset where $T$ is defined). We call $T$ {\it the $(C,F)$-transformation} associated with $(C_{n+1},F_n)_{n\ge 0}$. Below we will often use the following simple formulae
$$
\gathered
[A]_n\cap[B]_n=[A\cap B]_n,\quad
[A]_n\cup[B]_n=[A\cup B]_n,\\
[A]_n=\bigsqcup_{c\in C_{n+1}}[A+c]_{n+1},\quad
\mu([A+c]_{n+1})=\frac{\mu([A]_n)} {\# C_{n+1}},\\
T^r[A]_n=[r+A]_{n}
\endgathered
\tag1-1
$$
for all subsets $A,B\subset F_n$ and $r\in\Bbb Z$ such that $r+A\subset F_n$.

We note that in a similar way we can define  $(C,F)$-actions of an arbitrary countable discrete amenable group. In that case $(F_n)_{n\ge 0}$ must be a F{\o}lner sequence in the group.
The formulae \thetag{1-1} are all satisfied for arbitrary $(C,F)$-actions. In this paper we are mainly interested in $\Bbb Z$-actions. While proving Theorem~0.3 we will also need $(C,F)$-actions of $\Bbb Z^d$ which are Cartesian products of a fixed $\Bbb Z$-action.
A useful observation is that if $T$ is associated  with $(C_{n+1},F_n)_{n\ge 0}$ then the product  $\Bbb Z^d$-action
  $(T^{n_1}\times\cdots\times T^{n_d})_{(n_1,\dots,n_d)\in\Bbb Z^d}$ is associated with $(C_{n+1}^d,F_n^d)_{n\ge 0}$, where the upper index $d$ denotes the Cartesian power.

Another observation is that the $(C,F)$-construction is equivalent to the classical cut--and-stack construction of rank-one transformations.
Indeed, $X_n$ can be interpreted as the $n$-th tower consisting of the levels $[f]_n$, $f\in F_n$. It is cut into $\# C_{n+1}$ subtowers $[F_n+c]_{n+1}$, $c\in C_{n+1}$, which are then stack (with some spacers in-between) into a new, $(n+1)$-tower.
 $C_{n+1}$ is the set of locations of these subtowers inside the  $(n+1)$-tower.
More precisely, if we order the elements of $C_{n+1}$ as follows
$$
0=c(0)<\cdots < c(r_n-1)
 $$
then $c(i)$ is exactly the hight of the bottom level of the $(i+1)$-th subtower of $X_n$ inside $X_{n+1}$, $i=0,\dots, r_n-1$.
That is why in the following we will illustrate some abstract aspects of the $(C,F)$-construction with more common cut-and-stack pictures.

\head 2. High staircase construction
\endhead

In this section we prove Theorems~0.2 and 0.3.
We first formulate a couple of definitions.

\definition{Definition 2.1} Let $T$ be a measure preserving transformation of an infinite $\sigma$-finite measure space $(X,\goth B,\mu)$.
\roster
\item"\rom{(i)}"
A sequence of positive integers $a_n$ is called {\it mixing} for $T$ if
$\mu(T^{a_n}A\cap B)\to 0$ as $n\to \infty$ for all subsets $A,B\subset X$ of finite measure.
\item"\rom{(ii)}"
A sequence of intervals $[a_n,b_n)\subset\Bbb N$ is called {\it mixing} for $T$ if each sequence $d_n$ such that $a_n\le d_n<b_n$ is mixing for $T$.
\endroster
\enddefinition

Let $(z_n)_{n=1}^\infty$ and $(r_n)_{n=0}^\infty$ be two sequences of positive integers and $r_n\to\infty$ as $n\to\infty$.
We define $(C_n)_{n=1}^\infty$ and $(F_n)_{n=0}^\infty$ inductively by setting
$$
\align
& h_0\quad\text{is arbitrary},\\
& h_{n+1}:=r_n(h_n+z_n)+r_n(r_n-1)/2,\\
&c_{n+1}(0):=0, \quad c_{n+1}(i+1):=c_{n+1}(i)+h_n+z_n+i,\\
& C_{n+1}:=\{c_{n+1}(i)\mid i=0,\dots,r_n-1\},\\
& F_{n+1}:=\{0,\dots,h_{n+1}-1\}.
\endalign
$$

\definition{Definition 2.2}
The  $(C,F)$-transformation $T$ on a standard $\sigma$-finite measure space $(X,\goth B,\mu)$ associated with $(C_{n+1},F_n)_{n>0}$
will be called {\it a high staircase} (see Figure~2.1). If  $z_n=0$ for all $n$, we call $T$ a {\it pure staircase}.
\enddefinition

\midinsert
\beginpicture
\put {} at 0 0
\put {} at 0 215
\plot 0 0   0 190 20 190 20 195 40 195 40 200 60 200 60 205
80 205 80 210 100 210 100 215 120 215 120 0  0 0 /
\plot 0 140  120 140 /
\setdashes <3pt>
\plot 20 0 20 190 /
\plot 40 0 40 195 /
\plot 60 0 60 200 /
\plot 80 0 80 205 /
\plot 100 0 100 210 /
\setdots <2pt>
\plot 0 0 0 -15 /
\plot 120 0 120 -15 /
\setsolid
\arrow <5pt> [0.5, 1] from 0 -15 to 120 -15
\arrow <5pt> [0.5, 1] from 120 -15 to 0 -15
\put {$r_n$} at   60 -21
\setdots <2pt>
\plot 0 0   -15 0 /
\plot 0 190 -15 190 /
\plot 0 140 -15 140 /
\setsolid
\arrow <5pt> [0.5, 1] from -15 0 to -15 140
\arrow <5pt> [0.5, 1] from -15 140 to -15 0
\arrow <5pt> [0.5, 1] from -15 190 to -15 140
\arrow <5pt> [0.5, 1] from -15 140 to -15 190
\put {$h_n$} at -23 70
\put {$z_n$} at -23 165
\setdots <2pt>
\plot 120 0  135 0 /
\plot 0 190 135 190 /
\plot 120 215 135 215 /
\setsolid
\arrow <5pt> [0.5, 1] from 135 215 to 135 190
\arrow <5pt> [0.5, 1] from 135 190 to 135 215
\put {$r_n-1$} at   155 202
\arrow <5pt> [0.5, 1] from 135 190 to 135 0
\arrow <5pt> [0.5, 1] from 135 0 to 135 190
\put {$H_n$} at   145 95
\endpicture

\botcaption{Figure 2.1}
High staircase.
\endcaption
\endinsert

\definition{Definition 2.3}
By the  {\it restricted growth condition} we mean the following:
$$
\lim_{n\to\infty}\frac{r_n^2}{r_0r_1\cdots r_{n-1}}=0.
$$
\enddefinition

We note  that this condition implies
$$
\frac{r_n^2}{h_n}=\frac{r_n^2}{\mu(X_n)}\cdot \frac 1{r_0\cdots r_{n-1}}\to 0\quad\text{as }n\to\infty.\tag{2-1}
$$
Moreover, the restricted growth condition  is equivalent to  \thetag{2-1} in the case when $\mu$ is finite. We remark that Adams originally introduced it for the finite measure preserving pure staircases in the form of \thetag{2-1} \cite{Ad}.

We leave to the reader the proof of the following two simple lemmata.
The second one is an  $L^2$-version of \cite{Ad, Lemma~2.1}.

\proclaim{Lemma 2.4} Under the restricted growth condition a high staircase is infinite measure preserving if and only if \/ $\sum_{n=1}^\infty z_n/h_n=\infty$.
\endproclaim

\proclaim{Lemma 2.5} Given positive integers $R,L,r$ and a measurable set $B\subset X$ of finite measure, we have
$$
\bigg\|\frac 1R\sum_{i=0}^{R-1}U_T^{-i}1_B\bigg\|_2\le \bigg\|\frac 1L\sum_{i=0}^{L-1}U_T^{-ir}1_B\bigg\|_2 +\frac{rL}R\sqrt{\mu(B)}.
$$
\endproclaim

From now on we assume that $T$ is infinite measure preserving (see Lemma~2.4).
The following statement is a slight modification of  \cite{Ad, Lemma~2.2}.

\proclaim{Lemma 2.6} Let $[\alpha_n,\beta_n)$ be a sequence of intervals in $\Bbb Z$ which is mixing for each power of $T$ and $l_n\to\infty$. Take $k_n\in[\alpha_n,\beta_n)$. Then
$\frac 1{l_n}\sum_{i=0}^{l_n-1}U_T^{-k_ni}\to 0$ strongly as $n\to\infty$.
\endproclaim
\demo{Proof} Take a measurable subset $B\subset X$, $\mu(B)<1$. Then for each $l>0$,
$$
\bigg\|\frac 1l\sum_{i=0}^{l-1}U_T^{-ik_n}1_B\bigg\|_2^2=
\frac{\|1_B\|_2^2}{l}+\frac 1{l^2}\sum_{i\ne j=0}^{l-1}\mu(T^{(i-j)k_n}B\cap B).
$$
Since $\mu(T^{pk_n}B\cap B)\to 0$ for each $p\ne 0$ by the assumption of the lemma, there is $l$ such that
$$
\bigg\|\frac 1l\sum_{i=0}^{l-1}U_T^{-ik_n}1_B\bigg\|_2^2<\epsilon
$$
for all sufficiently large $n$.
 It remains to apply Lemma~2.5.
\qed
\enddemo

We now let $H_n:=h_n+z_n$.
The following lemma  is crucial in the proof of~Theorem~0.2.
\comment
It is worth to note that only the first claim  was used in \cite{Ad} to prove mixing of the probability preserving pure staircases.
For the high staircases it is not enough. That is why we add the second claim.
\endcomment

\proclaim{Lemma 2.7} If the restricted growth condition is satisfied
then
the sequence of intervals $[h_n,2H_n)$ is mixing for each non-zero power of $T$.
\endproclaim

\demo{Proof} Fix $j>0$.
We first show that  the sequence $(H_n)_{n>0}$ is mixing for $T^j$.
Take  subsets $A,B\subset F_n$.
Since
$$
jH_n+c_{n+1}(i)=c_{n+1}(i+j)-ji-\frac{j(j-1)}{2}\quad\text{if }i+j<r_n,
 $$
we have
$$
\align
\mu(T^{jH_n}[A]_n\cap[B]_n) &=\sum_{i=0}^{r_n-1}\mu(T^{jH_n} [A+c_{n+1}(i)]_{n+1}\cap[B]_n)\\
&=\sum_{i=0}^{r_n-j-1}\mu([-ij+A_j+c_{n+1}(i+j)]_{n+1}\cap[B]_n) +\bar o(1)\\
&=\frac 1{r_n}\sum_{i=0}^{r_n-j-1}\mu(T^{-ij}[A_j]_{n}\cap[B]_n)+\bar o(1),
\endalign
$$
where $A_j:=(A-{j(j-1)}/2)\cap F_n$. We note that $\bar o(1)$ here and below means ``uniformly small'' over all $A,B\subset F_n$ as $n\to\infty$.
Since $T$ is ergodic, $r_n^{-1}\sum_{i=0}^{r_n-1}U_T^{-i}\to 0$ strongly.
It follows that ${r_n}^{-1}\sum_{i=0}^{r_n-1}U_T^{-ij}\to 0$ strongly.
Hence we are done.

Now take $a_n\in [h_n,2H_n)$.
We are going to prove that $(a_n)_{n>0}$ is mixing for $T^j$.
Dropping to a subsequence we may assume  without loss of generality that $ja_n=kH_n+b_n$ for some  $0\le b_n<H_n$ and $0\le k\le 2j$.

Suppose first that $k\ne 0$.
Partition $A$ into three subsets $A_1, A_2,A_3$ such that $b_n+A_1\subset F_n$, $b_n+A_2\subset H_n+F_n$ and $b_n+A_3\subset [h_n,H_n)$.
We  verify mixing separately on each of these subsets.
We note first that $T^{b_n}[A_3]_n\cap [B]_n=\emptyset$,
$\mu([B]_n\triangle T^{-kH_n}[B]_n)=\bar o(1)$
and hence
$$
\mu(T^{ja_n}[A_3]_n\cap[B]_n)=\mu(T^{b_n}[A_3]_n\cap T^{-kH_n}[B]_n)=\bar o(1).
$$
To verify mixing  on on $A_1$ and $A_2$
it is enough to notice that
$$
\align
ja_n+A_1 &=kH_n+(b_n+A_1)\\
ja_n+A_2 &=(k+1)H_n+(b_n+A_2-H_n),
\endalign
$$
$(b_n+A_1)\cup(b_n+A_2-H_n)\subset F_n$ and use the fact that the sequence $(H_n)_{n>0}$ is mixing for both $T^k$ and $T^{k+1}$.

Now consider the second case when $k=0$.
Then $b_n\ge jh_n$.
Hence  we can partition $A$ into subsets $A_2$ and $A_3$ such that
$b_n+A_3\subset [h_n,H_n)$ and $b_n+A_2\subset H_n+F_n$.
Therefore $T^{ja_n}[A_3]_n\cap [B]_n=\emptyset$ and
$$
\mu(T^{ja_n}[A_2]_n\cap [B]_n)=\mu(T^{H_n}[b_n-H_n+A_2]_n\cap [B]_n)\to 0
$$
because $(H_n)_{n=1}^\infty$ is mixing for $T$.
\qed

\comment
Now let us show that  the sequence $(h_n)_{n>0}$ is mixing for $T^j$.
Passing, if necessary, to a subsequence we may assume without loss of generality that there are two integers $0\le k<j$ and $0\le b<h_n'$ such that $jz_n=kH_n+b_n$.
Then we have
$$
\align
jh_n+c_{n+1}(i)&=(j-k)H_n-b_n+c_{n+1}(i)\\
&=c_{n+1}(i+j-k)- (j-k)i-\frac{(j-k)(j-k-1)}2-b_n
\endalign
$$
if $i+j-k<r_n$.
Thus we come to what has been already considered above.

To show (i) take $a_n\in [h_n,2h_n)$ and prove that the sequence $(a_n)_{n>0}$ is mixing.  Dropping to a subsequence we may assume that
$$
ja_n=kh_n+b_n\quad\text{for some }0\le b_n<h_n\text{ and }j\le k\le 2j.
$$
Partition $A$ into two subsets $A_1, A_2$ such that $b_n+A_1\subset F_n$, $b_n+A_2\subset h_n+F_n$. Verify mixing separately for $A_1$ and $A_2$.
 For that it is enough to notice that
$$
\align
ja_n+A_1 &=kh_n+(b_n+A_1)\\
ja_n+A_2 &=(k+1)h_n+(b_n+A_2-h_n),
\endalign
$$
$(b_n+A_1)\cup(b_n+A_2-h_n)\subset F_n$ and use the fact that the sequence $(h_n)_{n>0}$ is mixing for both $T^k$ and $T^{k+1}$.
\endcomment

\enddemo

\demo{Proof of Theorem~0.2}
It is enough to prove that any sequence $(m_n)_{n=1}^\infty$ such that
$H_n\le m_n< H_{n+1}$ for all $n$ is mixing (or it contains a mixing subsequence).
We  find  integers $k_n$ and $t_n$ such that
$$
m_n=k_nH_n+t_n
$$
and $1\le k_n\le r_n$ and $0\le t_n<H_n$.
We set
$$
\align
C_{n+1}^1&:=\{c_{n+1}(i)\mid r_n-k_n\le i<r_n\},\\
C_{n+1}^2&:=\{c_{n+1}(i)\mid 0\le i<r_n-k_n\}.
\endalign
$$
{\bf Mixing on $F_n+C_{n+1}^1$.}
This corresponds to the domain $D_4$ on Picture~2.2.
Given $A,B\subset F_n$, we note that for each $0\le j<r_{n+1}-1$
$$
m_n+ A+C_{n+1}^1+c_{n+2}(j)= D-j+c_{n+2}(j+1),
$$
where $D:=t_n + A + C_{n+1}^1-H_{n+1}$.
It is important to notice that the subset $D-j$ is contained {\it essentially} in $F_{n+1}$, i.e. $\mu([(D-j)\cap F_{n+1}]_{n+1})=\mu([D]_{n+1})+\bar o(1)$. Hence
$$
\mu(T^{m_n}[A+C_{n+1}^1]_{n+1}\cap [B]_n)=\frac 1{r_{n+1}}\sum_{j=0}^{r_{n+1}-2}\mu(T^{-j}[D]_{n+1}\cap[B]_n)+\bar o(1)\to 0.
$$

\midinsert
\beginpicture
\put {} at 0 0
\put {} at 0 215
\plot 0 0   0 190 20 190 20 195 40 195 40 200 60 200 60 205
80 205 80 210 100 210 100 215 120 215 120 0  0 0 /
\plot 0 140  120 140 /
\setdots <2pt>
\plot 0 0 0 -15 /
\plot 0 0 0 -15 /
\plot 80 0 80 -15 /
\setsolid
\plot 80 0 80 140 /
\arrow <5pt> [0.5, 1] from 0 -15 to 80 -15
\arrow <5pt> [0.5, 1] from 80 -15 to 0 -15
\put {$r_n-k_n$} at   40 -21
\put {$D_4$} at 100 70
\put {$D_1$} at 40 30
\put {$D_3$} at 40 83
\put {$D_2$} at 40 123
\plot 0 60 80 60 /
\plot 0 110 80 110 /
\setdots <2pt>
\plot 0 0   -15 0 /
\plot 0 60   -15 60 /
\plot 0 110   -15 110 /
\plot 0 140 -15 140 /
\setsolid
\arrow <5pt> [0.5, 1] from -15 140 to -15 110
\arrow <5pt> [0.5, 1] from -15 110 to -15 140
\arrow <5pt> [0.5, 1] from -15 0 to -15 60
\arrow <5pt> [0.5, 1] from -15 110 to -15 60
\arrow <5pt> [0.5, 1] from -15 60 to -15 110
\arrow <5pt> [0.5, 1] from -15 140 to -15 0
\put {$h_n-t_n$} at -35 30
\put {$z_n$} at -25 85
\put {$t_n-z_n$} at -35 125
\endpicture
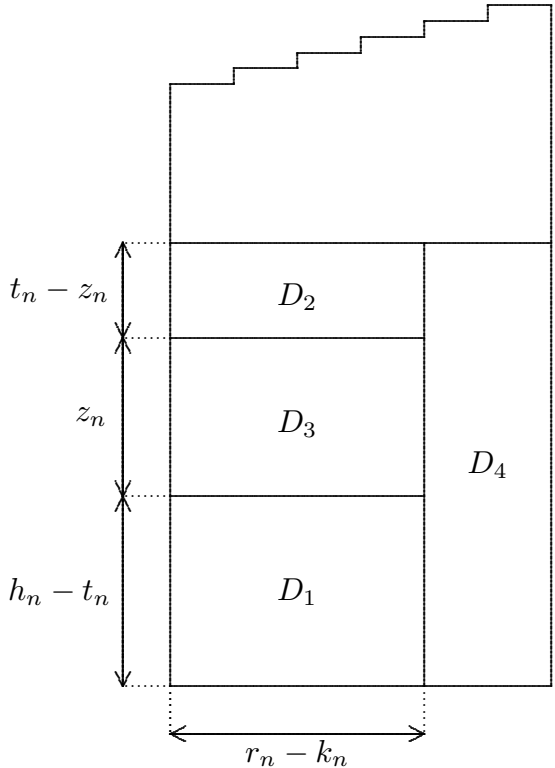
\botcaption{Figure 2.2} Partition of the $n$-tower into main domains.
\endcaption
\endinsert

\noindent{\bf Mixing on $F_n+C_{n+1}^2$.}
  Dropping to a further subsequence we can  assume that $k_n/r_n\to\delta_2$ for some $\delta_2<1$.
Indeed, if $k_n/r_n\to 1$ then
$$
\max_{D\subset F_n}\mu([D+C_{n+1}^2]_{n+1})/\mu([D]_n)\to 0
$$
and mixing on $F_n+C_{n+1}^2$ follows immediately.
Now take $c_{n+1}(i)\in C_{n+1}^2$.
We have
$$
m_n+c_{n+1}(i)=c_{n+1}(i+k_n)-k_ni-\frac{k_n(k_n-1)}2+t_n.
$$
Consider subsets $A,B\subset F_n$ such that $\mu([A\cup B]_n)<1$.
Partition $A$ into three subsets $A_1,A_2,A_3$ such that
$A_1+t_n\subset F_n$,  $A_2+t_n\subset H_n+F_n$  and $A_3+t_n\subset[h_n,H_n)$.
As for the graphical interpretation,  $A_i$ is the part of $A$ that lays inside the domain $D_i$, $1\le i\le 3$, on Figure~2.2.

As in the proof of Lemma~2.7, we verify mixing separately on every of these subsets.
However now $k_n$ may be unbounded and this fact essentially complicates the proof.
The restricted growth condition yields
$$
\max_{0\le i<r_n}\max_{D\subset F_n}|\mu([(D\pm k_ni\pm k_n(k_n-1)/2)\cap D]_n)-\mu([D]_n)|=\bar o(1),\tag2-2
$$
i.e. we may {\it neglect} rotations by $k_ni+k_n(k_n-1)/2$ inside $F_n$.
Taking this into account we  obtain mixing on $A_3$ in the same way as in the proof of Lemma~2.7.

Next, let us verify mixing on $A_1$. (Mixing on $A_2$ is verified in a similar way. We leave that to the reader.) Making use of \thetag{2-2} again, we obtain
$$
\mu(T^{m_n}[A_1+C_{n+1}^2]_{n+1}\cap [B]_n)
=\frac 1{r_n}\sum_{i=0}^{r_n-k_n-1}\mu(T^{-k_ni}[A_1']_n\cap [B]_n)+
\bar o(1),\tag2-3
$$
where $A_1':=(A_1+t_n- k_n(k_n-1)/2)\cap F_n$.

\midinsert

\beginpicture
\put {} at 0 0
\put {} at 0 215
\plot 0 0 120 0  120 160   0 160  0 0 /
\setplotsymbol({\llarge .})
\plot 1 0 79 0 /
\plot 40 140 59 140 /
\plot 60.5 120 79 120 /
\plot 80.5 100 99 100 /
\plot 100.5 80 119 80 /
\setplotsymbol({\tiny .})
\arrow <5pt> [0.5, 1] from 10 0 to 50 139
\arrow <5pt> [0.5, 1] from 30 0 to 70 119
\arrow <5pt> [0.5, 1] from 50 0 to 90 99
\arrow <5pt> [0.5, 1] from 70 0 to 110 79
\plot 0 160 0 190 20 190 20 195 40 195 40 200 60 200 60 205
80 205 80 210 100 210 100 215 120 215 120 160 /
\setdots <2pt>
\plot 0 0 0 -15 /
\plot 80 0 80 -15 /
\setsolid
\arrow <5pt> [0.5, 1] from 0 -15 to 80 -15
\arrow <5pt> [0.5, 1] from 80 -15 to 0 -15
\put {$r_n-k_n$} at   40 -21
\setdots <2pt>
\plot 60 140  135 140 /
\plot 120 80 135 80 /
\setsolid
\arrow <5pt> [0.5, 1] from 135 140 to 135 80
\arrow <5pt> [0.5, 1] from 135 80 to 135 140
\put {$(r_n-k_n)k_n$} at   163 110
\setdots <2pt>
\plot 40 140  -15 140 /
\plot 60 120 -15 120 /
\setsolid
\arrow <5pt> [0.5, 1] from -15 140 to -15 120
\arrow <5pt> [0.5, 1] from -15 120 to -15 140
\put {$k_n$} at -23 131
\setdashes <3pt>
\plot 20 0 20 190 /
\plot 40 0 40 195 /
\plot 60 0 60 200 /
\plot 80 0 80 205 /
\plot 100 0 100 210 /
\endpicture
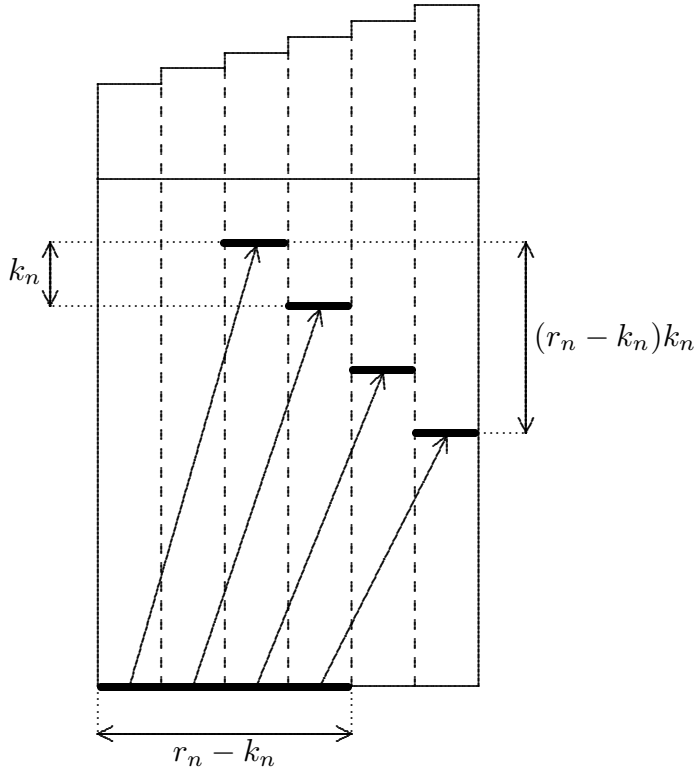
\botcaption{Figure 2.3} Image of the bottom layer of $D_1$ under $T^{m_n}$.
\endcaption
\endinsert

There exists an integer $p_n$ such that
$$
h_{p_n}\le (r_n-k_n)k_n<h_{p_{n}+1}.
$$
Since $\delta_2\ne 1$ and $r_n\to\infty$, it follows that  $p_n\to\infty$.
Consider now separately three possible cases.

{\bf (A)} Let $\frac{(r_n-k_n)k_n}{h_{p_n}}\to\infty$ and $k_n\le h_{p_n}$ for all $n$.
Then we can find $k_n'$ such that $h_{p_n}\le k_nk_n'<2h_{p_n}$.
It follows that
$\frac{r_n-k_n}{k_n'}\to\infty$.
We can choose a sequence $(l_n)_{n=1}^\infty$ so that
$l_n\to\infty$ and $\frac{r_n-k_n}{l_nk_n'}\to\infty$.
By Lemmata~2.6 and 2.7,
$$
\frac 1{l_n}\sum_{i=0}^{l_n-1}U_T^{-ik_nk_n'}\to 0\quad\text{strongly as }n\to\infty.
$$
Now we deduce from Lemma~2.5 that
$$
\frac 1{r_n}\sum_{i=0}^{r_n-k_n-1}\mu(T^{-k_ni}[A_1']_n\cap [B]_n)\le
\bigg\|\frac {1}{l_n}\sum_{i=0}^{l_n-1}U_T^{-k_nk_n'i}1_{[B]_n}\bigg\|_2
+\frac{l_nk_n'}{r_n-k_n}.
$$
Mixing on $A_1$ in this case follows from the above inequality and \thetag{2-3}.

Before we pass to the remaining two cases let us show that
  $$
{H_n^2}/{h_{n+1}}\to\infty.\tag2-4
$$
 Indeed,
$$
\frac{H_nr_n+r_n(r_n-1)/2}{H_n^2}
\le \frac{r_n}{H_n}+\frac{r_n^2}
{H_n^2}\le \frac{r_n}{h_n}+\frac{r_n^2}{h_n^2}\to 0
$$
in view of \thetag{2-1}.

{\bf (B)} Let $\frac{(r_n-k_n)k_n}{h_{p_n}}\to\infty$ and $k_n> h_{p_n}$ for all $n$.
If $H_{p_n}<k_n$ for all sufficiently large $n$ then
$$
\frac{2\delta_2}{1-\delta_2}>\frac{k_n^2}{(r_n-k_n)k_n}> \frac{H_{p_n}^2}{(r_n-k_n)k_n}>\frac{H_{p_n}^2}{h_{p_n+1}}\to\infty
$$
according to \thetag{2-4}, a contradiction.
Hence, passing to a subsequence we may assume that
$$
h_{p_n}\le k_n\le H_{p_n}
$$
for all sufficiently large $n$.
Then by Lemmata~2.7 and 2.6 we deduce from \thetag{2-3} that
$\mu(T^{m_n}[A_1+C_{n+1}^2]_{n+1}\cap [B]_n)\to 0$.

{\bf (C)} Consider now the remaining third case where $\frac{(r_n-k_n)k_n}{h_{p_n}}$ does not go to infinity. Then  dropping to a  subsequence we may assume that there is  $\delta_1\ge 1$ with
$$
\frac{(r_n-k_n)k_n}{h_{p_n}}\to \delta_1\quad \text{as }n\to\infty.
$$
 It follows that there is $\delta>0$ such that
$$
h_{p_n}>\delta k_n^2\quad \text{for all sufficiently large $n$}.\tag2-5
$$
We deduce from \thetag{2-4} that
 $H_{p_n-1}^2/{h_{p_n}}\to\infty$. This plus \thetag{2-5} imply
$H_{p_n-1}/k_n\to\infty$. Hence $H_{p_n-1}>k_n>h_{p_n-1}$ for all sufficiently large $n$.
Therefore again we can find an integer $k_n^*$ such that $H_{p_n-1}\le k_n^*k_n< 2H_{p_n-1}$.
The sequence of intervals $[H_{p_n},2H_{p_n})$ is mixing for all powers of $T$ by Lemma~2.7
We also have
$$
\frac{r_n-k_n}{k_n^*}=\frac{(r_n-k_n)k_n}{k_n^*k_n}\ge
\frac{(r_n-k_n)k_n}{2H_{p_n-1}}\ge \frac{\delta_1}4\cdot\frac{h_{p_n}}{H_{p_n-1}}\to\infty.
$$
Therefore arguing as in {\bf (A)} we obtain mixing in the case {\bf (C)}.
\qed
\enddemo

\remark{Remark \rom{2.8}} We note that in the particular case when $z_n/h_n\to 0$ as $n\to\infty$ the proof of Theorem~0.1 is simplified significantly.
It can be carried out as a slight modification of the reasoning of Adams in \cite{Ad}.
In particular, we then no longer need to consider the cases (B) and (C).
Furthermore, we can ``reduce'' the interval $[h_n,2H_n)$ in the statement of Lemma~2.7 to the interval $[h_n,2h_n)$.
Nevertheless,  this particular case of Theorem~0.2 is enough to demonstrate Theorems~0.3 and 0.4. However in the proof of our main result---Theorem~0.1---we essentially use Theorem~0.2 in the full generality. More precisely, it is utilized to satisfy \thetag{3-2} below.
\endremark

Our next purpose is to construct a power weakly mixing rank-one infinite measure preserving transformation which is mixing. For that we will use high staircases and Theorem~0.2.

\demo{Proof of Theorem 0.3}
We will construct these high staircases via an inductive procedure.
Suppose that after   $p$  steps  we have already defined $F_0,C_1,F_1,\dots, C_{n_p}, F_{n_p}$.

{\it Step $p+1$}. Consider an auxiliary finite measure preserving pure staircase  $T_{p+1}$ associated with a sequence $(C_{k}^{(p+1)},F_{k}^{(p+1)})_{k\ge 0}$
such that $F_0^{(p+1)}:=F_{n_p}$.
Assume that the restricted growth condition is satisfied for $T_{p+1}$.
Let $(X^{(p+1)},\mu_{p+1})$ be the space of this action.
We normalize $\mu_{p+1}$ in such a way that
$$
\mu_{p+1}([0]_0)=\frac 1{\# C_1\cdots\# C_{n_p}}.
$$
Take any finite sequence $l=(l_1,\dots,l_{p+1})$ of non-negative integers such that
$\|l\|:=\max_{0<i\le p+1}|l_i|\le p+1$.
Since $T_{p+1}$ is  mixing \cite{Ad}, the  transformation
$$
S_l:=T_{p+1}^{l_1}\times\cdots\times T_{p+1}^{l_{p+1}}
$$
of the product space $(X^{(p+1)},\mu_{p+1})^{p+1}$ is ergodic.
Then there are $N_{p+1}>0$ and $M_{p+1}>0$ such that for all disjoint subsets $A,B\subset (F_0^{(p+1)})^{p+1}$ of equal cardinality there exist subsets $A'\subset[A]_0$ and $B'\subset[B]_0$ and their partitions $A'=\bigsqcup_{i=1}^{M_{p+1}}A_i$, $B'=\bigsqcup_{i=1}^{M_{p+1}}B_i$ such that
$$
\gathered
\mu_{p+1}^{\times(p+1)}(A')>0.5\mu_{p+1}^{\times{p+1}}([A]_0),\quad
\mu_{p+1}^{\times{p+1}}(B')>0.5\mu_{p+1}^{\times{p+1}}([B]_0),\\
\text{$A_i,B_i$ are $N_{p+1}$-cylinders (some may be empty) and}\\
S_l^iA_i=B_i\quad\text{for all }i=1,\dots,M_{p+1}.
\endgathered
\tag2-6
$$
It is assumed that $N_{p+1}$ and $M_{p+1}$ are common for all  $S_l$ with $\|l\|\le p+1$.
We now set
$$
C_{n_p+1}:=C^{(p+1)}_1, F_{n_p+1}:=F^{(p+1)}_1,\dots, F_{n_p+N_{p+1}}:=F^{(p+1)}_{N_{p+1}}.
$$
Next, we add several more subsets $(C_{i},F_{i})_{i=n_p+N_{p+1}+1}^{n_{p+1}}$ in the high staircase shape and such that the corresponding parameters $z_i$ incorporated into $C_i$ are {\it large}.
These subsets are needed to get infinite measure {\it in the limit}.
The $(p+1)$-step is now completed.

Continuing this procedure  infinitely many times, we obtain the entire sequence $(C_{i+1},F_{i})_{i=0}^\infty$.
Denote by $(X,\mu,T)$ the associated $(C,F)$-dynamical system.
We normalize $\mu$ in such a way that $\mu([0]_0)=1$.
By construction, $T$ is a high staircase.
Since $z_i$ are {\it large} for infinitely many $i$, it follows that $\mu(X)=\infty$.
Of course, we may assume without loss of generality that the restricted growth condition is satisfied.
Hence $T$ is mixing by Theorem~0.2.

To verify that $T$ is power weakly mixing it is enough to notice that
$$
\gather
F_0^{(p+1)}=F_{n_p}, \ F_{N_{p+1}}^{(p+1)}=F_{n_p+N_{p+1}}\quad\text{and}
\\
\mu_{p+1}([A]_{N_{p+1}})=\mu([A]_{n_p+N_{p+1}})\tag2-7
\endgather
$$
for each  subset $A\subset F_{N_{p+1}}^{(p+1)}$, $p>0$,  and use \thetag{2-6}.
We note that \thetag{2-7} follows from the normalization conditions for $\mu_{p+1}$ and $\mu$.
\qed

\enddemo

\head 3.
Spectral multiplicities of mixing infinite measure preserving transformations
\endhead

In this section we use Theorem~0.2 to show Theorems~0.1 and 0.4.
We first prove an auxiliary lemma about about cyclic spaces for products of unitary operators.

\proclaim{Lemma 3.1}
Let $U, V$ be two unitary operators in a Hilbert space $\Cal H$
 and let $V$ has a  simple spectrum.
Let $\Cal C$ be a $U$-cyclic space generated by a vector $h_U$ and let $h_V$ be a cyclic vector for $V$.
 If there is  a sequence $n_i\to\infty$ such that $U^{n_i}\to I$ and $V^{n_i}\to aV^*$ weakly for some $a>0$ then $\Cal C\otimes \Cal H$ is the $U\otimes V$-cyclic space generated by $h_U\otimes h_V$.
\endproclaim
\demo{Proof}
Let $h_U$ and $h_V$ be cyclic vectors for $U$ and $V$ respectively.
Denote by $\Cal D$ the $(U\otimes V)$-cyclic space generated by $h_U\otimes h_V$.
Since $(U\otimes V)^{n_i}\to aI\otimes V^*$ and $\Cal D$ is invariant under the weak limits of powers of $U\otimes V$, it follows that $h_U\otimes V^*h_V\in\Cal D$.
In a similar way,  $h_U\otimes V^nh_V\in\Cal D$ for each $n\in\Bbb Z$. Therefore $h_U\otimes\Cal H\subset\Cal D$ which implies, in turn, that $\Cal C\otimes\Cal H\subset\Cal D$.
The converse inclusion is obvious.
\qed
\enddemo

We need to consider a class of transformations which is  more general than the high staircases.
Let three sequences of positive integers $(z_n)_{n=1}^\infty$, $(r_n)_{n=1}^\infty$ and
$(d_n)_{n=1}^\infty$ be given such that $d_n\le r_n$ and $\lim_{n\to\infty}d_n/r_n=\delta\ge 0$.
Suppose that the following are satisfied
$$
\align
& C_{n+1}:=\{c_{n+1}(i)\mid i=0,\dots,r_n-1\},\\
&0=c_{n+1}(0)<c_{n+1}(1)<\cdots<c_{n+1}(r_n-1), \\ &c_{n+1}(i+1):=c_{n+1}(i)+h_n+z_n+i-d_n \quad\text{if }i=d_n,\dots,r_n-1,
\endalign
$$
We call the transformation associated with  $(C_{n+1}, F_n)_{n\ge 0}$   a {\it $(1-\delta)$-partially high staircase}.
Geometrically this  means that on the $n$-th step we cut the $n$-th tower into $r_n$ subtowers and arrange spacers on the tops of $d_n$ first subtowers in an arbitrary way. On the tops of the remaining $r_n-d_n$ subtowers the spacers are arranged in the high staircase way.
If $\delta=0$ then we call the transformation an {\it almost high staircase}.
Since the proportion of the first $d_n$ subtowers of such a transformation goes to  $0$,
it can be deduced from the proof of Theorem~0.2 that every almost high staircase satisfying the restricted growth condition is mixing.

\demo{Proof of Theorem 0.1}
Suppose first that $\infty\not\in M$.
It was shown in \cite{DaR} that  there exists an ergodic conservative infinite measure preserving transformation $S$ such that $\Cal M(S)=M$
and $U_S^{i_m}\to I$ weakly for  a certain sequence $i_m\to\infty$.
Let $(Y,\nu)$ denote the $\sigma$-finite space where $S$ acts.

Now we will ``force mixing'' of such transformations along the scheme outlined in the introduction.
 Fix  a sequence $\delta_n\to 0$.
We will construct a sequence of $(1-\delta_n)$-partially high staircases $T_n$ associated with $(C_k^{(n)},F_{k-1}^{(n)})_{k>0}$ such that
$$
U_{T_n}^{H_k^{(n)}}\to\delta_n U_{T_n}^*\quad\text{weakly as }k\to\infty,
\tag3-1
$$
 where $(H_k^{(n)})_{k=1}^\infty$ is a subsequence  of $(i_m)_{m=1}^\infty$.
For that we put
$$
\align
C_{k+1}^{(n)}&=\{c^{(n)}_{k+1}(i)\mid 0\le i<r_k^{(n)}\},\quad\text{where}\\
c^{(n)}_{k+1}(i+1)&:=
\cases
0, & \text{for }i=-1\\
c^{(n)}_{k+1}(i)+H_{k}^{(n)}+1, &\text{for } 0\le i< d_k^{(n)}-1\\
c^{(n)}_{k+1}(i)+H_{k}^{(n)}+i-d_k^{(n)}+1, &\text{for } d_k^{(n)}-1\le i< r_k^{(n)}-1
\endcases
\endalign
$$
and $H_{k}^{(n)} :=h_k^{(n)}+z_k^{(n)}$ and $F_{k+1}^{(n)}=[0,h_{k+1}^{(n)})$, where
$$
h_{k+1}^{(n)} :=r_k^{(n)}H_{k+1}^{(n)}+d_k^{(n)}-1+ \frac{(r_k^{(n)}-d_k^{(n)}-1)(r_k^{(n)}-d_k^{(n)})}{2}.
$$
The initial {\it heights}  $h_0^{(n)}$ (i.e. sets $F_0^{(n)}$) are not specified yet.
This will be done below.
We now impose   the following restrictions on the parameters $r_k^{(n)},d_k^{(n)},z_k^{(n)}$:
$$
\gather
r_k^{(n)}\to\infty,\quad2\delta_n\ge{d_k^{(n)}}/{r_k^{(n)}}\to\delta_n \quad\text{as }k\to\infty,\quad\\
 \sum_{k=1}^\infty z_k^{(n)}/h_k^{(n)}=\infty\quad \text{and}\\
\{H_k^{(n)}\mid k\in\Bbb N\}\subset \{i_m\mid m\in\Bbb N\}.\tag3-2
\endgather
$$
Several additional restrictions will appear below.
To verify \thetag{3-1} take subsets $A,B\subset [r_k^{(n)},h_k^{(n)})$.
Then up to $\bar o(1)$ in measure the subset $T_n^{H_k^{(n)}}[A]_k$ equals
$$
\bigsqcup_{i=1}^{d_k^{(n)}-1} [A-1+c_{k+1}^{(n)}(i)]_{k+1}\sqcup
\bigsqcup_{i=1}^{r_k^{(n)}-d_k^{(n)}-1} [A-i+c_{k+1}^{(n)}(i+d_k^{(n)})]_{k+1}.
$$
Hence
$$
\langle U_{T_n}^{H_k^{(n)}}1_{[A]_k},1_{[B]_k}\rangle
=\frac{d_k^{(n)}}{r_k^{(n)}}
\langle U_{T_n}^*1_{[A]_k},1_{[B]_k}\rangle+
\frac{1}{r_k^{(n)}}\sum_{i=1}^{r_k^{(n)}-d_k^{(n)}-1}\langle U_{T_n}^{-i}1_{[A]_k},1_{[B]_k}\rangle+\bar o(1)
$$
and \thetag{3-1} follows because $T_n$ is ergodic.

The sought-for transformation $T$ will appear as a $(C,F)$-transformation associated with a {\it concatenated} sequence
$$
F_0^{(1)},C_1^{(1)},\dots,F_{k_1}^{(1)}, C_1^{(2)}, F_1^{(2)}, C_2^{(2)},\dots,F_{k_2}^{(2)}, \dots.\tag3-3
 $$
Of course,  to make this concatenation well defined we need to satisfy a {\it compatibility} condition
$F_{k_n}^{(n)}=F_0^{(n+1)}$.
By this condition we determine all the initial heights $h_0^{(n)}$ except for $h_0^{(1)}$ which is chosen in an arbitrary way.
It remains to specify the sequence of {\it stopping times} $(k_n)_{n=1}^\infty$.
This will be done inductively.

Fix a vector $v$ in $L^2(Y,\nu)$ and denote by $\Cal C$ the $U_S$-cyclic subspace generated by $v$.
Fix a dense countable subset $(v_i)_{i=1}^\infty$ in $\Cal C$.
Suppose that we have already determined $k_1,\dots,k_{n-1}$.
Then $F_0^{(n)}$ is defined by the compatibility condition.
The other sets $C_1^{(n)},F_1^{(n)},C_2^{(n)},\dots$ are defined by the above recurrent formulae.
Let $T_n$ be the associated $(C,F)$-transformation acting on a space
$(X^{(n)},\mu^{(n)})$.
As in the proof of Theorem~0.3 we normalize $\mu^{(n)}$ in such a way that
$$
\mu^{(n)}([0]_0)= \bigg({\prod_{i=1}^{n-1}\#C^{(i)}_1\cdots\#C^{(i)}_{k_i}}\bigg)^{-1}.
$$
Since $T_n$ has a simple spectrum, we can choose a cyclic vector
$w_n$  for $U_{T_n}$.
 Since \thetag{3-1} and \thetag{3-2} are satisfied,  $v\otimes w_n$ is a cyclic vector for $(U_S\restriction\Cal C)\otimes U_{T_n}$ by Lemma~3.1.
Hence we can select a large $k_n$ in such a way that there are a
subset $A\subset F_{k_n}^{(n)}$, a function  $w_n'$ in the linear span of the {\it indicators} $1_{[f]_{k_n}}$, $f\in A$, and $M>0$ such that $M+A\subset F_{k_n}^{(n)}$ and
$$
\max_{f'\in F_0^{(n)}}
\Bigg\|v_i\otimes 1_{[f']_0}-\sum_{j=-M}^M\alpha_{i,j}(U_S\otimes U_{T_n})^j\,v\otimes w_n'\Bigg\|_2<\epsilon_n \tag3-4
$$
for all $i=1,\dots,n$, where $\alpha_{i,j}$ are some real numbers and $\epsilon_n\to 0$ very fast.
This completes the $n$-th step of the inductive procedure.
On this step we defined a fragment $C_1^{(n)},\dots,F^{(n)}_{k_n}$
of \thetag{3-3}.
Continuing this infinitely many times we construct the entire sequence \thetag{3-3}.
Rename it as $F_0,C_1,F_1,\dots$.
Let $T$ be the associated $(C,F)$-transformation  and let $(X,\mu)$ be the space of this transformation.
The measure $\mu$ is normalized so that $\mu([0]_0)=1$.
By the construction, $T$ is an almost high  staircase and $\mu(X)=\infty$.
By choosing the parameters $r_{k}^{(n)}$ in a right way we may assume without loss of generality  that the restricted growth condition is satisfied for $T$.
Hence $T$ is mixing.

Now we are going to show that $\Cal M(S\times T)=\Cal M(S)$.
Let $\Cal H^{(n)}$ denote the linear span of the indicators $1_{[f]_0}$, $f\in F_0^{(n)}$
in $L^2(X^{(n)},\mu^{(n)})$ and let
$\Cal H_{k}$ denote the linear span of the indicators $1_{[f]_k}$, $f\in F_k$,
in $L^2(X,\mu)$.
We note that $F_{k_1+\cdots+k_n}=F_0^{(n)}$ for all $n>0$.
Moreover, it is easy to see that the natural identification
$$
1_{[f]_{k_1+\cdots+k_n}}\leftrightarrow 1_{[f]_0},\quad f\in F_{k_1+\cdots+k_n},
$$
extends to  a linear isomorphism of $\Cal H_{k_1+\cdots+k_n}$ onto ${\Cal H}^{(n)}$.
Moreover, it is an isometry because of the normalizations imposed on
 $\mu$ and $\mu^{(n)}$.
Hence \thetag{3-4} yields  the existence of $w_n'\in\Cal H_{k_1+\cdots+k_{n+1}}$ such that
$$
\max_{f'\in F_{k_1+\dots+k_n}}
\Bigg\|v_i\otimes 1_{[f']_{k_1+\dots+k_n}}-\sum_{j=-M}^M\alpha_{i,j}(U_S\otimes U_T)^jv\otimes w_n'\Bigg\|_2<\epsilon_n
$$
for all $i=1,\dots,n$, $n\in\Bbb N$.
It follows that the subspace $\Cal C\otimes L^2(X,\mu)$
 is cyclic for $U_S\otimes U_T$.

Consider now a decomposition of $L^2(Y,\nu)$ into $U_S$-cyclic spaces
$\Cal C_{i,j}$:
$$
L^2(Y,\nu)=\bigoplus_{i\in\Cal M(S)}\bigoplus_{j=1}^i\Cal C_{i,j}
$$
such that $\sigma_{i,j}\sim \sigma_{i,j'}$ for all $i,j,j'$ and $\sigma_{i,j}\perp \sigma_{i',j'}$ if $i\ne i'$, where $\sigma_{i,j}$
is a measure of maximal spectral type of $U_S\restriction\Cal C_{i,j}$.
It follows that a subspace $C_{i,j}\oplus C_{i',j'}$ is $U_S$-cyclic if $i\ne i'$.
Given a $U_S$-cyclic space $\Cal C$, we showed above how to construct a mixing infinite measure preserving almost high staircase dynamical system $(X,\mu,T)$ such that $\Cal C\otimes L^2(X,\mu)$ is $U_S\otimes U_T$-cyclic.
It is clear that the construction can be obviously modified so to obtain the same for an arbitrary countable family of $U_S$-cyclic subspaces.
In particular, we can construct
a mixing infinite measure preserving almost high staircase system $(X,\mu,T)$ such that the subspaces
$$
(\Cal C_{i,j}\oplus\Cal C_{i',j'})\otimes L^2(X,\mu)\quad\text{are $U_S\otimes U_T$-cyclic}\tag3-5
$$
for all  $i\ne i'\in\Cal M(S)$, $1\le j\le i$ and $1\le j\le i'$.
Denote by $\widehat\sigma_{i,j}$ a measure of maximal spectral type for
$(U_S\restriction\Cal C_{i,j})\otimes U_T$.
It follows from \thetag{3-5} that $\widehat\sigma_{i,j}\perp \widehat\sigma_{i',j'}$ if $i\ne i'$.
Of course, $\sigma_{i,j}\sim \sigma_{i,j'}$ for all $i,j,j'$.
Since
$$
L^2(Y,\nu)\otimes L^2(X,\mu)=\bigoplus_{i\in\Cal M(S)}\bigoplus_{j=1}^i\Cal C_{i,j}\otimes L^2(X,\mu)
$$
and $\Cal C_{i,j}\otimes L^2(X,\mu)$ is $U_S\otimes U_T$-cyclic for each pair $i,j$, it follows that $\Cal M(S\times T)=\Cal M(S)$.

Thus the theorem is ``almost'' proved. It  only remains unclear  whether $T\times S$ is ergodic or not.
Therefore to obtain the desired ergodicity we will modify the construction of $T$: we will force ergodicity for  $T\times S$ simultaneously with forcing mixing for $T$.
For that we alternate the above argument with the argument from the proof of Theorem~0.3.
On the odd steps we construct fragments of \thetag{3-3} exactly as above
to retain the desired spectral properties and mixing.
It remains to explain the construction on the even steps.

Suppose we have already constructed
$$
F_0^{(1)}, C_1^{(1)},\dots,
C^{(n-1)}_{k_{n-1}},F^{(n-1)}_{k_{n-1}}\tag{3-6}
$$
with $n-1$ odd.
Now we consider a pure staircase infinite $(C,F)$-sequence
$$
F_0^{(n)}, C_1^{(n)},F_1^{(n)}, C_2^{(n)},\dots
$$
such that $F_0^{(n)}=F_{k_n-1}^{n-1}$.
Recall that this means that
$$
C^{(n)}_k=\{c_k(i)\mid i=0,\dots,r_k-1\},
$$
$c_k(0)=0$, $c_k(i+1)=c_k(i)+h_{k-1}+i$, where $h_{k-1}$ is the hight of $F_{k-1}^{(n)}$ and $h_k=r_kh_{k-1}+r_k(r_k-1)/2$.
We also assume that the restricted growth condition $r_k^2/h_k\to 0$ as $k\to\infty$ is satisfied.
Let $T_n$ be the $(C,F)$-transformation acting on a measured space $(X^{(n)},\mu^{(n)})$.
Notice that $\mu^{(n)}(X^{(n)})<\infty$.
As above, we normalize $\mu_n$ in such a way that
$$
\mu^{(n)}([0]_0)= \bigg({\prod_{i=1}^{n-1}\#C^{(i)}_1\cdots\#C^{(i)}_{k_i}}\bigg)^{-1}.
$$
By \cite{Ad}, $T_n$ is mixing.
Hence the Cartesian product $S\times T_n$ is ergodic \cite{FW}.
Recall that $S$ is a $(C,F)$-map \cite{DaR}.
Let $(C_{m+1}^S,F_m^S)_{m=0}^\infty$ denote the corresponding sequence.
Since $S\times T_n$ is ergodic, there are integers $k_n>0$ and $M_n>0$ such that for all disjoint subsets $A,B\subset F_0^S\times F_0^{(n)}$ of equal cardinality there exist subsets $A'\subset[A]_0$ and $B'\subset[B]_0$ and their partitions $A'=\bigsqcup_{i=1}^{M_n}A_i$ and
$B'=\bigsqcup_{i=1}^{M_n}B_i$ such that
$$
\gathered
\nu\times\mu^{(n)}(A')>   0.5(\nu\times\mu^{(n)})([A]_0),\\
\nu\times\mu^{(n)}(B')>   0.5(\nu\times\mu^{(n)})([B]_0),\\
A_i,B_i\text{ are $k_n$-cylinders (some may be empty) and}\\
(S\times T_n)^i A_i=B_i\text{ for all }i=1,\dots,M_n.
\endgathered
\tag3-7
$$
Then we ``continue'' \thetag{3-6} with the following fragment:
$C_1^{(n)}, F_1^{(n)},\dots,C_{k_n}^{(n)}, F_{k_n}^{(n)}$.

Thus we explained the construction procedure entirely.
The corresponding $(C,F)$-transformation $T$  is as desired.
Indeed, we have $\Cal M(S\times T)=\Cal M(S)$ and $S\times T$ is mixing due to the odd steps of our construction process.
Furthermore, \thetag{3-7} implies that $S\times T$ is ergodic.

Now let $\infty\in M$.
  Assume first that the set $M':=M\setminus\{\infty\}$ is nonempty.
This case would immediately come to the previous one whenever we know that the maximal spectral type of $S\times T$ is  singular.
Unfortunately, we do not see easy ways to show this property for every $S$.
However, it holds for certain $S$ that `splits into product' of two other transformations.
Thus our goal is to show that $S$ can always be chosen in this special way.
For that we proceed in several steps.

Fix an ergodic conservative  rigid transformation $S$  with $\Cal M(S)=M'$ \cite{DaR}.
We also fix a family of $S$-cyclic spaces $\Cal C_{i,j}$, $i\in M'$, $1\le j\le i$ as above.

{\it Claim A.} There exists an infinite measure preserving transformation $R$ of a $\sigma$-finite standard measure space $(Z,\goth F,\kappa)$ such that
\roster
\item"$(A1)$"
$R$ is rigid along a subsequence of $(i_m)_{m>0}$ (the latter is a rigidity sequence for $S$),
\item"$(A2)$"
$R\times S$ is conservative and ergodic,
\item"$(A3)$"
the subspaces
$L^2(Z,\kappa)\otimes (\Cal C_{i,j}\oplus\Cal C_{i',j'})$
are $U_R\otimes U_S$-cyclic
for all  $i\ne i'\in\Cal M(S)$, $1\le j\le i$ and $1\le j\le i'$.
\endroster

Recall that the group Aut$(Z,\kappa)$ of all $\mu$-preserving invertible transformations of $(Z,\kappa)$ is Polish in the weak topology \cite{CK}.
It is well known (and easy to verify) that the set $\Cal T_1$ of all transformations $R$ satisfying $(A1)$ is a dense $G_\delta$ in Aut$(Z,\kappa)$.
The map
$$
f_{i,i',j,j'}:R\mapsto U_R\otimes (U_S\restriction (\Cal C_{i,j}\oplus\Cal C_{i',j'}))
$$
from
Aut$(Z,\kappa)$ to  the group of unitary operators (of a corresponding Hilbert space) equipped with the weak operator topology is continuous for each quadruple of indices $i,i',j,j'$.
Recall that the set of unitary operators with a simple continuous spectrum is a  $G_\delta$ in the the group of all unitary operators \cite{Na}.
Hence the set $\Cal T_3$ of all transformations $R$ satisfying $(A3)$ is a $G_\delta$ in Aut$(Z,\kappa)$.
Since the set of ergodic conservative transformations is a $G_\delta$ in Aut$(Z\times Y,\kappa\times\nu)$ \cite{CK} and the map $R\mapsto R\times S$ is continuous, it follows that the set $\Cal T_2$ of all transformations $R$ satisfying~$(A2)$ is also a $G_\delta$ in Aut$(Z,\kappa)$.
As follows from the first part of the proof of Theorem~0.1, the intersection $\Cal T_2\cap\Cal T_3$ is nonempty (the transformation $T$ constructed there belongs to it).
Since the intersection is invariant under the conjugacy and the conjugacy class of each conservative ergodic transformation is dense in Aut$(Z,\kappa)$ \cite{CK}, we deduce that $\Cal T_2\cap\Cal T_3$ is a dense $G_\delta$ in Aut$(Z,\kappa)$.
Hence so is $\Cal T_1\cap\Cal T_2\cap\Cal T_3$.
The claim follows.

We deduce from Claim A that the product transformation $R\times S$ is ergodic, conservative and rigid along a subsequence of $(i_m)_{m>0}$.
Moreover,
$$
L^2(Z\times Y,\kappa\times\nu)=\bigoplus_{i\in M'}\bigoplus_{1\le j\le i} L^2(Z,\kappa)\otimes \Cal C_{i,j},
$$
and the restrictions of $U_{R\times S}$ to  cyclic subspaces $L^2(Z,\kappa)\otimes \Cal C_{i,j}$ and
$L^2(Z,\kappa)\otimes \Cal C_{i',j'}$ are  either unitarily equivalent if $i=i'$ or spectrally disjoint if $i\ne i'$.
In particular,
$\Cal M(R\times S)=\Cal M(S)=M'$.

{\it Claim B.} There is a mixing transformation $T$ such that
\roster
\item"$(B1)$"
$
L^2(Z\times Y\times X,\kappa\times\nu\times\mu)=\bigoplus_{i\in M'}\bigoplus_{1\le j\le i} L^2(Z,\kappa)\otimes \Cal C_{i,j}\otimes L^2(X,\mu),
$
\item"$(B2)$"
every subspace $\Cal C_{i,j}':=L^2(Z,\kappa)\otimes \Cal C_{i,j}\otimes L^2(X,\mu)$ is
$U_{R\times S\times T}$-cyclic,
\item"$(B3)$" the unitary operators $U_{R\times S\times T}\restriction \Cal C_{i,j}'$ and $U_{R\times S\times T}\restriction \Cal C_{i',j'}'$ are either unitarily equivalent if $i=i'$ or spectrally disjoint if $i\ne i'$.
\item"$(B4)$" $R\times S\times T$ is ergodic and conservative.
\endroster
This claim can be shown by repeating almost verbally the argument from the first part of the proof of Theorem~0.1 (for $M$ without $\infty$).
One should just replace $S$ with $R\times S$.

The following assertion is well known. We state it without proof.

{\it Claim C.} If $V$ and $W$ are unitary operators  acting in infinite dimensional  Hilbert spaces $\Cal H_V$  and $\Cal H_W$ such that $\infty\not\in\Cal M(V\otimes W)$ then the maximal spectral type of $V$ (and $W$) is singular.

{\it Final step.} We deduce from $(B2)$ and $(B3)$ that $\Cal C_{i,j}'':=\Cal C_{i,j}\otimes L^2(X,\mu)$ is a $U_{S\times T}$-cyclic subspace for each pair $i,j$.
Furthermore, the unitary operators $U_{S\times T}\restriction \Cal C_{i,j}''$  and $U_{S\times T}\restriction \Cal C_{i',j'}''$ are unitarily equivalent if $i=i'$ or spectrally disjoint if $i\ne i'$.
Hence $\Cal M(S\times T)=M'$.
Next, it follows from Claim B that $\Cal M(R\times S\times T)=M'$.
Then Claim~C yields that the maximal spectral type of $U_{S\times T}$ is singular.
Therefore if $B$ is a (probability preserving) Bernoulli shift then
$$
\Cal M(S\times T\times B)=\Cal M(S\times T)\cup\{\infty\}=M.
$$
It follows from $(B4)$ and \cite{FW} that
the transformation $S\times T\times B$ is ergodic.
Since the measure on the space of  this transformation is non-atomic,
the transformation is conservative.

It remains to consider the case where $M=\{\infty\}$.
However  this case is well known.
We recall that an ergodic conservative infinite measure preserving $K$-transformation
has Lebesgue spectrum of infinite multiplicity \cite{Pa}.
Any such
a transformation is mixing \cite{KS}. 
\qed
\enddemo

Our next purpose is to prove Theorem~0.4.
For that we need an auxiliary lemma from \cite{Ag} and \cite{Ry2}.

\proclaim{Lemma 3.2}
Let $V$  be unitary operator with a simple spectrum in Hilbert space $\Cal H$. Assume  that for each $n\in\Bbb N$, there are a sequence $k_t^{(n)}\to\infty$ and    reals  $\alpha_n,\beta_n>0$ such that
$V^{k_t^{(n)}}\to \alpha_nI+\beta_nV^*$ weakly as $t\to\infty$,
 and $\alpha_1,\alpha_2,\dots$ are pairwise different.
Then $\exp(V)$ has a simple spectrum.
\endproclaim

We note that $\exp(V)$ has a simple spectrum if and only if each symmetric power $V^{\odot n}$ has a simple spectrum and the convolution powers of a measure of maximal spectral types for $V$ are all mutually singular.

\demo{Proof of Theorem 0.4}
We only sketch the idea of the proof.
As in the proof of Theorem~0.1 we will use the forcing mixing technique.
For that construct inductively a sequence of $(1-\alpha_n-\beta_n)$-partially high staircases $T_n$ such that
the weak closure of the powers of $U_{T_n}$ contains $\alpha_nI+\beta_nU_{T_n}^*$ with $\alpha_n,\beta_n \to 0$ as $n\to\infty$.
By Lemma~3.2, all the symmetric powers $U_{T_n}^{\odot k}$ have a simple spectrum,
$k=1,2,\dots$.
Let $T_n$ be associated with a sequence
$(C^{(n)}_{k},F^{(n)}_{k-1})_{k\ge 1}$.
Then we construct a $(C,F)$-transformation $T$ associated to
$$
F_0^{(1)},C_1^{(1)},\dots,F_{k_1}^{(1)}, C_1^{(2)}, F_1^{(2)}\dots,F_{k_2}^{(2)}, \dots,
 $$
where the stopping times $k_1,k_2,\dots$ are chosen in such a way to retain the simple spectrum of the symmetric powers ``in the limit'', i.e. for $T$. This is possible because the property to have a simple spectrum is  approximable: if a unitary operator admits a sequence of cyclic spaces {\it appoximating} the entire space then the operator has a simple spectrum.
We note that $T$ is an almost high staircase since $\alpha_n,\beta_n\to 0$.
 Hence it is conservative and ergodic.
 Moreover, we can satisfy the restricted growth condition.
Hence $T$ is mixing by Theorem~0.2.
\qed
\enddemo

\head 4. Applications to Poisson suspensions
\endhead

Let $(X,\goth B)$ be a standard Borel space and let $\mu$ be an infinite $\sigma$-finite non-atomic measure on $\goth B$.
Fix an increasing sequence of Borel subsets $X_1\subset X_2\subset\cdots$
with $\bigcup_{i>0}X_i=X$ and $\mu(X_i)<\infty$ for each $i$.
A Borel subset is called {\it bounded} if it is contained in some $X_i$.
Let $\widetilde X_i$ denote the space of finite measures on $X_i$.
For each bounded subset $A\subset X_i$, let $N_A$ stand for the map
$$
\widetilde X_i\ni\omega\mapsto \omega(A)\in\Bbb R.
$$
Denote by  $\widetilde {\goth B}_i$ the smallest $\sigma$-algebra on $\widetilde X_i$ in which all the maps  $N_A$, $A\in\goth B\cap X_i$, are measurable.
It is well known that $(\widetilde X_i,\widetilde {\goth B}_i)$
is a standard Borel space.
Denote by $(\widetilde X,\widetilde{\goth B})$ the projective limit of the sequence
$$
(\widetilde X_1, \widetilde{\goth B}_1)\leftarrow(\widetilde X_2,\widetilde{\goth B}_2)\leftarrow\cdots,
 $$
where the arrows denote the (Borel) natural {\it restriction} maps.
Then $(\widetilde X,\widetilde{\goth B})$ is a standard Borel space.
To put it in other way, $\widetilde X$ is the space of measures on $X$ which are $\sigma$-finite  along $(X_i)_{i>0}$.
We define a probability measure $\widetilde\mu$ on  $(\widetilde X,\widetilde{\goth B})$ by the following two conditions:
\roster
\item"\rom{(i)}"
$\widetilde\mu\circ N_A^{-1}=e^{-\mu(A)}\sum_{n=0}^\infty\frac{\mu(A)^n}{n!}\delta_n$ for each bounded subset $A$, where $\delta_n$ is the  Dirac measure on $\Bbb R$
supported at $n$,
\item"\rom{(ii)}"
if $A$ and $B$ are disjoint bounded subsets of $X$ then the random variables $N_A$ and $N_B$ on $(\widetilde X,\widetilde{\goth B},\widetilde\mu)$ are independent.
\endroster
These conditions determine $\widetilde\mu$ in a unique way.
 Let $T$ be a Borel  invertible transformation $T$  such that $T$ and $T^{-1}$ preserve the subclass of bounded subsets.
If $T$ preserves $\mu$ then it induces a Borel isomorphism $\widetilde T$ of $\widetilde X$ by the formula $\widetilde T\widetilde \omega:=\widetilde \omega\circ T$.
We recall that the dynamical system $(\widetilde X,\widetilde{\goth B},\widetilde\mu,\widetilde T)$ is called the {\it Poisson suspension} of $(X,\goth B,\mu,T)$ \cite{CFS}.
A probability preserving transformation is called a {\it Poissonian automorphism} if it is isomorphic to the  Poissonian suspension of an infinite measure preserving transformation.
If $T$ has no invariant subsets of positive and finite measure then $\widetilde T$ is weakly mixing.

We  state some new spectral realization problems:
\roster
\item"{\bf (P1)}" which spectral multiplicities\footnote{In the finite measure preserving case we consider the spectral multiplicities of operators restricted to the orthocomplement to the constants.} are realizable on weakly mixing Poissonian
automorphisms?
\item"{\bf (P2)}"
which spectral multiplicities are realizable on mixing Poissonian
automorphisms?
\endroster

The same questions can be also asked for the Gaussian automorphisms.
Since each Poissonian automorphism is spectrally equivalent to some Gaussian one, any (partial) solution of (P1) or (P2) is also a solution in the class of Gaussian automorphisms.
It is unclear whether the converse holds.

  It is known that $\Cal M(T)$ is either $\{1\}$ or infinite for each Gaussian (and therefore Poissonian) automorphism \cite{CFS}.
 Recall that the first mixing Gaussian automorphism with a simple spectrum appeared in \cite{New}.

Since $U_{\widetilde T}=\exp(U_T)$ \cite{Ne}, Theorem~0.4 implies the existence of a mixing Poissonian automorphism $\widetilde T$ with a simple spectrum (Corollary~0.5).
It is worth to note that $\widetilde T$ possesses some other interesting properties: if an invertible transformation $S$ preserves $\widetilde\mu$ and commutes with $\widetilde T$ then $S=\widetilde R$ for a $\mu$-preserving transformation $R$ commuting with $T$ \cite{Ro2}.
Also, $\widetilde T$ does not split into a Cartesian product of two other transformations \cite{Ro2}.
We note that the existence of non-mixing Poissonian automorphisms with a simple spectrum was mentioned in \cite{Ro2}.

\example{Example 4.2}
We will show that for each $p>1$, the semigroup $\{p,p^2, p^3,\dots\}$
is realizable as the set of spectral multiplicities for a
mixing Poissonian automorphism.
Indeed, let $I_p$ stand for the identity
transformation on $\Bbb Z/p\Bbb Z$.
Let $T$ be a mixing infinite measure preserving transformation $T$
 such that $\exp(T)$ has a simple spectrum (see Theorem~0.4).
The Cartesian product $T\times I_p$  is no longer ergodic.
However it  is mixing.
Hence the Poisson suspension $\widetilde{T\times I_p}$ is also mixing.
It follows from \cite{Ry3} that  $\Cal M(\widetilde{T\times I_p})=\{p,p^2,p^3,\dots\}$.
 \endexample

Infinite ``non-semigroups'' of positive integers also can appear
as Poissonian multiplicities.

\example{Example 4.3} Let $T$ be as in Example~4.2.
Then the transformation $T\odot T$ is mixing and
therefore its Poissonian suspension $\widetilde{T\odot T}$
is also mixing.
Moreover,
$$
U_{\widetilde{T\odot T}}=\exp(U_{T\odot T})=\exp(U_T\odot U_T)=\bigoplus_{n=0}^\infty
(U_T\odot U_T)^{\odot n}.
$$
Since $\Cal M((U_T\odot U_T)^{\odot n})=\{(2n)!/(2^nn!)\}$ and
the measure of spectral types of $(U_T\odot U_T)^{\odot n}$ are pairwise
disjoint, we obtain that
$$
\Cal M(\widetilde{T\odot T})=\bigg\{\frac{(2n)!}{2^nn!}\,\bigg|\,n\in\Bbb N\bigg\}=\{1,3,3\cdot 5,3\cdot 5\cdot 7,\dots\}.
$$
\endexample

Instead of $\widetilde{T\odot T}$ one can also consider Poissonian suspensions of other {\it natural  factors}  $T^{\otimes k}/ \Gamma$ of the product $T^{\otimes k}$, where $\Gamma$ is a subgroup of the full symmetric group $S_k$.
We leave details to the reader.

\head 5. Concluding remarks
\endhead

{\bf 5.1.} There exist other mixing infinite measure preserving constructions which can be used in this paper instead of high staircases:
\roster
\item"$(\circ)$" {\it High stochastic constructions with vanishing deterministic part}.
 No principal difficulties arise to adapt Ornstein's random spacer techniques from \cite{Or}   to the infinite setup.
\item"$(\circ)$" {\it Pure staircases and almost pure staircases with rapidly growing sequence of cuts}.
It follows from some unpublished results of the second named author that
the restricted growth condition can be waived from the statement of Theorem~0.2 provided that $r_n/h_n\to 0$ and $r_n\to\infty$.
If, in addition, $\sum_{n=1}^\infty r_n/h_n=\infty$ then the corresponding $(C,F)$-transformations are infinite measure preserving.
In this case we can restrict ourself to the pure  (and almost pure) staircases.
\endroster
However, not striving for the full generality, we choose the  high staircases satisfying the restricted growth condition as the most effective and fast way leading to the proof of Theorems~0.1--0.4.

{\bf 5.2.} We note that the {\it stopping times} appearing in the proofs of Theorems~0.1, 0.3 and 0.4 when we force mixing or ergodicity are not defined effectively. The same takes place in the other applications of the method to force mixing in \cite{Ag}, \cite{Ry2}, \cite{Da4}.
We hope to improve the method in a future work by developing an explicit effective algorithm of selecting the stopping times.

{\bf 5.3.} This remark is for aesthetically minded readers who may do not like the appearance of  {\it almost} high staircases in the proof of Theorem~0.1. They will be pleased to know that the almost high staircases there can be replaced  completely with high staircases (for the expense of some complication of the   proof, of course).
We now briefly outline the proof of this claim.
Let $S$ be as in the proof of Theorem~0.1.
In particular, $S^{i_m}\to I$ as $m\to\infty$.
We  will construct high staircase transformations $T_n$ such that
for for some $L_n$,
$$
\gather
U_{T_n}^{H^{(n)}_{q,k}}\to \Cal P_q(U_{T_n}):=\frac 1 {q+1}(I+U_{T_n}^{-1}+\cdots+U_{T_n}^{-q}) \quad\text{as $k\to\infty$ and}
\tag 5-1\\
\text{$(H^{(n)}_{q,k})_{k=1}^\infty$  is a subsequence of $(i_m)_{m=1}^\infty$ }\tag5-2
\endgather
$$
for each $q\ge L_n$.
To satisfy \thetag{5-1}, every number $q\ge L_n$ must occur in the sequence $(r^{(n)}_k)_{k=1}^\infty$ infinitely many times.
The property \thetag{5-2} is satisfied by choosing $(z^{(n)}_{k})_{k>0}$ in a right way.
Now let $v_S$ be a $U_S$-cyclic vector in a cyclic subspace $\Cal C\subset L^2(Y,\nu)$.
Let $v_{T_n}$ be a $U_{T_n}$-cyclic vector in $L^2(X_n,\mu_n)$.
Denote by $\widetilde {\Cal C}$ the $U_S\otimes U_{T_n}$-cyclic space generated by $v_S\otimes v_{T_n}$.
It follows from \thetag{5-1} and \thetag{5-2} that
$$
v_S\otimes \Cal P_q(U_{T_n})v_{T_n}\in\widetilde{\Cal C}\quad\text{for all $q\ge L_n$}.
$$
Hence
$$
(q+1)v_S\otimes \Cal P_{q}(U_{T_n})v_{T_n}-
qv_S\otimes \Cal P_{q-1}(U_{T_n})v_{T_n}=v_S\otimes U_{T_n}^{-q}v_{T_n}
\in\widetilde {\Cal C}
$$
for each $q> L_n$.
In particular, given any $p\in\Bbb Z$, we have
$v_S\otimes U_{T_n}^{p-i_{m}}v_{T_n}\in\widetilde {\Cal C}$ for all sufficiently large $m$.
Hence $U_S^{i_m}v_S\otimes U_{T_n}^pv_{T_n}\in\widetilde {\Cal C}$.
Passing to the (weak) limit as $m\to\infty$, we obtain that $v_S\otimes U_{T_n}^pv_{T_n}\in\widetilde {\Cal C}$.
This yields
$$
\widetilde {\Cal C}=\Cal C\otimes L^2(X_n,\mu_n),\tag5-3
$$ as desired.
Now to force mixing we must retain \thetag{5-3} in the limit as $n\to\infty$ for a {\it representative} family of $U_S$-cyclic subspaces $\Cal C$. The sought-for  high staircase $T$ appears as a certain limit of the sequence $(T_n)_{n>0}$.
At the same time we must have $L_n\to\infty$ in order to obtain $r_n\to\infty$ for $T$.
The rest of the argument is as in the proof of Theorem~0.1.

\enddocument